\documentclass[]{svproc}
\usepackage{url}

\usepackage{mathtools}
\usepackage{amsfonts}
\usepackage{bbm}
\usepackage{amsmath}
\usepackage{amssymb}
\DeclarePairedDelimiterX{\inp}[2]{\langle}{\rangle}{#1, #2}
\usepackage{xcolor}

\newtheorem{assumption}{Assumption}

\begin{document}
\mainmatter              
%
\title{Bridging Koopman Operator and time-series auto-correlation based Hilbert-Schmidt operator}
\titlerunning{Koopman and time-series}  
%
\author{Yicun Zhen\inst{1} \and Bertrand Chapron\inst{1} \and Etienne M\'{e}min\inst{2}}
\authorrunning{Zhen, Chapron and M\'{e}min} 
%
\tocauthor{Yicun Zhen, Bertrand Chapron, and Etienne M\'{e}min}
\institute{Institut Français de Recherche pour l'Exploitation de la Mer, Plouzan\'{e} 29280, France\\
\email{zhenyicun@protonmail.com}
\and
INRIA/IRMAR Campus universitaire de Beaulieu, 35042 Rennes Cedex, France. 
}

\maketitle              

\begin{abstract}
Given a stationary continuous-time process $f(t)$,  the Hilbert-Schmidt operator $A_{\tau}$ can be defined for every finite $\tau$\cite{Vautard1989SingularSA}. Let $\lambda_{\tau,i}$ be the eigenvalues of $A_{\tau}$ with descending order. In this article, a Hilbert space $\mathcal{H}_f$ and the (time-shift) continuous one-parameter semigroup of isometries $\mathcal{K}^s$ are defined. Let $\{v_i, i\in\mathbb{N}\}$ be the eigenvectors of $\mathcal{K}^s$ for all $s\geq 0$. Let $f = \displaystyle\sum_{i=1}^{\infty}a_iv_i + f^{\perp}$ be the orthogonal decomposition with descending $|a_i|$. We prove that $\displaystyle\lim_{\tau\to\infty}\lambda_{\tau,i} = |a_i|^2$. The continuous one-parameter semigroup $\{\mathcal{K}^s: s\geq 0\}$ is equivalent, almost surely, to the classical Koopman one-parameter semigroup defined on $L^2(X,\nu)$, if the dynamical system is ergodic and has invariant measure $\nu$ on the phase space $X$. 
\end{abstract}

\section{Introduction}
Let $\{f(t)\in\mathbb{C}: t\geq 0\}$ be a continuous time process. We assume that $f$ has zero temporal mean and the lagged moments exist for all $s\geq 0$:
\begin{align}
    \rho(s) := \lim_{T\to\infty}\frac{1}{T}\int_{0}^{T}f(t)\bar{f}(t+s)\mathtt{d}t. \label{eq: def rho s}
\end{align}
Define $\rho_{-s} = \bar{\rho}_{s}$.
In \cite{Vautard1989SingularSA} the self-adjoint operator $A_{\tau}$ is defined to act on $L^2([0,\tau])$:
\begin{align}
    (A_{\tau}g)(t) = \frac{1}{\tau}\int_{0}^{\tau} g(s)\rho(t-s) \mathtt{d}s, \label{eq: A tau}
\end{align}
for every $g\in L^2([0,\tau])$, and for all $t\in [0,\tau]$. When $\rho\in L^2_{\text{loc}}(\mathbb{R})$ and $\rho(s)\neq 0$ for almost all $s\in [0,\tau]$ , $A_{\tau}$ is a Hilbert-Schmidt operator. In particular, $A_{\tau}$ is compact and always has a purely  punctual spectrum. As stated in \cite{Vautard1989SingularSA}, the singular spectrum analysis (SSA) algorithm is based on the spectral analysis of $A_{\tau}$.

While in practice the SSA method has been applied successfully to a large variety of time series, in a theoretical purpose, yet with practical consequences, one may ask ourselves what is the relation between $A_{\tau_1}$ and $A_{\tau_2}$ for different $\tau_1$ and $\tau_2$? And what is the asymptotic behavior of $A_{\tau}$ as $\tau\to\infty$? In what way is the spectral property of $A_{\tau}$ related to  intrinsic properties of the dynamical system? In this article we generalize the idea and tools developed in \cite{Zhen2021} and apply them to study of $A_{\tau}$. We shall prove that
\begin{align}
    \lim_{\tau\to\infty}\lambda_{\tau,i} = |a_i|^2,
\end{align}
where $\lambda_{\tau,i}$ is the $i-$th largest eigenvalue of $A_{\tau}$ and $a_i$ is the $i-$th largest (in modulus) coefficient of some eigenvector $v_i$ (of unit length) of the time-shift operator $\mathcal{K}^s$ (for all $s\geq 0$) in the orthogonal decomposition of $f$:
\begin{align}
    f = \sum_{i=1}^{\infty}a_iv_i + f^{\perp},\label{eq: f decompose}
\end{align}
where $f^{\perp}$ denotes the the expression of $f$ in the orthogonal complement of the space spanned by the time-shift operator eigenfunctions.
If there are only finitely many $i$ (say only $N$ terms in the summation) in Eq.\eqref{eq: f decompose}, then we set $a_i = 0$ for $i>N$. The time-shift operator $\mathcal{K}^s$ is closely related to the classical Koopman operator, which is defined to act, as a time-shift operator, on some function space whose domain is the whole phase space of the dynamical system.

In section 2 we present the main result and a brief introduction of the mathematical tools used by the proof of the main result. All the quantities mentioned above are defined rigorously in section 2. The detailed proof of the main result is presented in section 3.

\paragraph{Notes and Comments.}
The main result as well as the techniques and ideas used for the proof  are close in spirit to those developed in \cite{Zhen2021}. 
However, the Hilbert-Schmidt operator $A_{\tau}$ is defined for continuous time process and the theory developed in \cite{Zhen2021} does not cover the continuous-time case. The objective of this paper is to confirm that the asymptotic behavior of the Hilbert-Schmidt operator $A_{\tau}$ is well related to Koopman theory.
\section{Preliminaries and the main result}

Let $\{f(t): t\geq 0\}$ be a continuous-time process. 
\begin{assumption}\label{assumption: stationary}
Assume that 
\begin{align}
    \lim_{T\to\infty}\frac{1}{T}\int_{0}^Tf(t)\mathtt{d}t = 0,
\end{align}
and that $\rho(s)$ is well-defined by Eq.\eqref{eq: def rho s} for all $s\geq 0$.
\end{assumption}
For any $s\geq 0$, we use $F_s$ to denote the time series $\{F_s(t) = f(t+s): t\geq 0\}$. For any two time series $g = \{g(t):t\geq 0\}$ and $h = \{h(t): t\geq 0\}$, we define the new time series 
\begin{align}
    ag+bh = \{ag(t) + bh(t): t\geq 0\},
\end{align}
where $a,b\in\mathbb{C}$. We consider the following linear space:
\begin{align}
  \widetilde{\mathcal{H}}_f = \text{Span}_{\mathbb{C}}\{F_{s}: s\geq 0\}.
\end{align}
Each element  $h\in\widetilde{\mathcal{H}}_s$ can be written as 
\begin{align}
    h = \sum_{i=1}^n c_iF_{s_i},
\end{align}
for any $n\geq 1, c_i\in\mathbb{C}, s_i\geq 0$. The existence of $\rho(s)$ allows us to define the following positive semi-definite Hermitian form:
\begin{align}
    \inp{h}{g} = \lim_{T\to\infty}\frac{1}{T}\int_{0}^Th(t)\bar{g}(t)\mathtt{d}t.
\end{align}
Let $V = \{v\in\widetilde{\mathcal{H}}_f: \inp{v}{v} = 0\}$. Since the Hermitian form is positive semi-definite, $V$ is a linear subspace of $\tilde{\mathcal{H}}_f$. And the Hermitian form is strictly positive-definite on the quotient space $\widetilde{\mathcal{H}}_f/V$. Hence it defines an inner product on $\widetilde{H}_f/V$. We define 
\begin{align}
    \mathcal{H}_f := \overline{\widetilde{\mathcal{H}}_f/V}
\end{align}
where the closure is taken with respect to the inner product defined above.

We define the operator $\mathcal{K}^s$ on $\widetilde{\mathcal{H}}_f$ for any $s, s_1\geq 0$:
\begin{align}
    \mathcal{K}^s F_{s_1} = F_{s_1 + s}.
\end{align}
It is obvious that 
\begin{align}
    \inp{\mathcal{K}^sh}{\mathcal{K}^sg} = \inp{h}{g},
\end{align}
for any $h,g\in\tilde{\mathcal{H}}_f$ and any $s\geq 0$. Hence $\mathcal{K}^s$ is well-defined on $\widetilde{\mathcal{H}}_f/V$, and can be further extended to the whole $\mathcal{H}_f$ by continuity. Therefore we obtain a one parameter family of isometric operators $\mathcal{K}^s$ that acts on the Hilbert space $\mathcal{H}_f$. And obviously we have 
\begin{align}
    \mathcal{K}^{s_1}\mathcal{K}^{s_2} = \mathcal{K}^{s_1+s_2}.
\end{align}
To simplify the notation, we use $f$ to also denote the continuous-time process $F_0$.  We further assume that
\begin{assumption}\label{assumption: K continuous}
\begin{align}
\lim_{s\to 0^+}\|\mathcal{K}^s f - f\|_{\mathcal{H}_f} = 0.
\end{align}
\end{assumption}
In other words, assumption \ref{assumption: K continuous} assumes that the curve: 
\begin{align}
    \gamma: [0,\infty)&\rightarrow \mathcal{H}_f\nonumber \\
     t&\rightarrow \mathcal{K}^tf
\end{align} is continuous. Since $\mathcal{H}_f$ is generated by $f$ and $\mathcal{K}^s$ are isometries for all $s\geq 0$, assumption \ref{assumption: K continuous} implies that $\mathcal{K}^s\to I$ in the strong operator topology as $s\to 0^+$. In other words, $\{\mathcal{K}^s: s\geq 0\}$  forms a strongly continuous semigroup of isometries on $\mathcal{H}_f$.

Under assumption \ref{assumption: K continuous}, we have the following decomposition theorem (see Theorem 9.3 in \cite{SzkefalviNagy1970HarmonicAO} ).
\begin{theorem}\label{theorem: wold continuous time}
 Let $\{\mathcal{K}^s: s\geq 0\}$ be a strongly continuous  semigroup of isometries on a Hilbert space $\mathcal{H}$. Then $\mathcal{H}$ has the orthogonal decomposition $\mathcal{H} = \mathcal{H}_U\bigoplus\mathcal{H}_{NU}$, where $\mathcal{H}_U = \displaystyle\bigcap_{s\geq 0}\mathcal{K}^s\mathcal{H}$, and $\mathcal{H}_{NU}$ is isomorphic to $L^2([0,\infty], \mathcal{H}_{0})$ for some Hilbert space $\mathcal{H}_0$. $\mathcal{H}_U$ and $\mathcal{H}_{NU}$ are invariant under $\mathcal{K}^s$ for all $s\geq 0$. The operator $\mathcal{K}^s$ restricted on $\mathcal{H}_U$ is a strongly continuous semigroup of unitary operators. And $\mathcal{K}^s$ restricted to $\mathcal{H}_{NU}$ acts as the unilateral shift operator, i.e. for any $\gamma\in \mathcal{H}_{NU} = L^2([0,\infty], \mathcal{H}_0)$, 
\begin{align}
    (\mathcal{K}^s\gamma)(t) = \gamma(t+s)\in\mathcal{H}_0.
\end{align}
\end{theorem}
Theorem \ref{theorem: wold continuous time}  provides us with an useful tool to deal with the completely nonunitary component of $\mathcal{K}^s$. For the unitary component, we have the following spectral representation theorem.
\begin{theorem}\label{theorem: spectral representation continuous}
Let $\{U(s): s\geq 0\}$ be a strongly continuous semigroup of unitary operators on a Hilbert space $\mathcal{H}$. Assume that $\mathcal{H}$ can be generated by $U$ and some $f\in\mathcal{H}$. Then there exists a unitary map $\phi: \mathcal{H}\rightarrow L^2(\mathbb{R}, d\mu)$ where $\mu$ is some positive finite measure on $\mathbb{R}$, such that 
\begin{align}
    (\phi(f))(x) &= 1,\\
    (\phi(\mathcal{K}^sg))(x) &= e^{isx}(\phi(g))(x)
\end{align}
for all $g\in\mathcal{H}$, $x\in\mathbb{R}$, and $s\geq 0$.
\end{theorem}
 Theorems \ref{theorem: wold continuous time} and \ref{theorem: spectral representation continuous} suggest the orthogonal decomposition $\mathcal{H}_f = \mathcal{H}_{f,U}\bigoplus\mathcal{H}_{f,NU} = L^2(\mathbb{R},d\mu_f)\bigoplus L^2([0,\infty], \mathcal{H}_{f,0})$. Furthermore, we can write $\mu_f = \mu_{f,d} + \mu_{f,c}$, where $\mu_{f,d}$ is a countable sum of Dirac measures and $\mu_{f,c}$ is continuous with respect to the Lebesgue measure. $\mu_{f,c}$ can be composed both of an absolutely continuous part and a singular continuous part. The decomposition of $\mu_f$ suggests the orthogonal decomposition $\mathcal{H}_{f,U} = L^2(\mathbb{R},d\mu_{f,d})\bigoplus L^2(\mathbb{R},d\mu_{f,c})$. In sum, we have
 \begin{align}
     f = f_{NU} + f_{d} + f_c \label{eq: f = fd + fc + fNU},
 \end{align}
where $f_{NU}\in L^2([0,\infty], \mathcal{H}_{f,0})$, $f_d\in L^2(\mathbb{R}, d\mu_{f,d})$, and $f_{c}\in L^2(\mathbb{R},d\mu_{f,c})$. Note that these subspaces are pair-wise orthogonal and are all invariant under $\mathcal{K}^s$ for all $s\geq 0$. The support of $\mu_{f,d}$ consists of countably many points. Each point $x_i$ in the support of $\mu_{f,d}$ corresponds to an eigenvector $v_i\in\mathcal{H}_f$ of $\mathcal{K}^s$ for all $s\geq 0$, i.e.
\begin{align}
    (\phi(a_iv_i))(x) = \begin{cases}
    1 &\text{if $x = x_i$,}\\
    0 &\text{otherwise,}
    \end{cases}\label{eq: phi(vi) = dirac}
\end{align}
and $\mu_{f,d}(\{x_i\}) = |a_i|^2$, where $a_i$'s are the coefficients of the eigenvectors in the following decomposition: 
\begin{align}
    f = \sum_{i}a_iv_i + f_{NU} + f_{c}.\label{eq: f = aivi + fNU + fc}
\end{align}
We rearrange the index of $v_i$ so that $|a_1|\geq |a_2|\geq \cdots \geq 0$. In order to make connection with $A_{\tau}$, we need the following lemmas.

\begin{lemma}\label{lemma: integral exists}
For any $\tau> 0$ and any $g\in L^2([0,\tau])$, the following integral
\begin{align}
    \int_{0}^{\tau}g(s)\mathcal{K}^sf \mathtt{d}s
\end{align}
is well-defined and is an element of $\mathcal{H}_f$.
\end{lemma}
The proof of this and the following lemma use standard argument from mathematical analysis and we leave the proof to the interested readers.

Let 
\begin{align}
    \widetilde{\mathcal{H}}_f^{\text{int}} = \{\int_{0}^{\tau}g(s)\mathcal{K}^sf\mathtt{d}s: \tau > 0, g\in L^2([0,\tau])\}.
\end{align}
$\widetilde{\mathcal{H}}^{\text{int}}_f$ is a linear subspace of $\mathcal{H}_f$. We have 
\begin{lemma}
\begin{align}
\overline{\widetilde{\mathcal{H}}^{\text{int}}_f} = \mathcal{H}_f.    
\end{align}
\end{lemma}


For simplicity, we use the notation $L^2_{\tau} := L^2([0,\tau])$. Given lemma \ref{lemma: integral exists}, for any $g_1, g_2\in L^2([0,\tau])$ and $t\in [0,\tau]$, 
we define the Hermitian form ${\bf A}_{\tau}: L^2_{\tau}\times L^2_\tau\rightarrow \mathbb{C}$:
\begin{align}
    {\bf A}_{\tau}(g_1,g_2) = \frac{1}{\tau}\inp[\Big]{\int_{0}^{\tau}g_1(t)\mathcal{K}^tf\mathtt{d}t}{\int_{0}^{\tau}g_2(s)\mathcal{K}^sf\mathtt{d}s}_{\mathcal{H}_f}.\label{eq: <gAg> = <int,int>}
\end{align}
Cauchy-Schwartz inequality implies that
\begin{align}
    |{\bf A}_\tau (g_1,g_2)|^2&\leq \frac{1}{\tau^2}{\Big \|}\int_{0}^\tau g_1(s)\mathcal{K}^sf\mathtt{d}s{\Big \|}_{\mathcal{H}_f}^2{\Big \|}\int_0^\tau g_2(s)\mathcal{K}^s f\mathtt{d}s{\Big \|}^2\\
    &\leq \frac{1}{\tau^2} \|g_1\|_{L^2_\tau}^2\|g_2\|_{L^2_\tau}^2\|f\|_{\mathcal{H}_f}^4,
\end{align}
where $\inp{}{}_{L^2_\tau}$ refers to the inner product in $L^2_{\tau}$ and $\inp{}{}_{\mathcal{H}_f}$ refers to the inner product in $\mathcal{H}_f$. Therefore Riesz representation theorem warrants that there exists a linear bounded operator $A_\tau: L^2_\tau\rightarrow L^2_\tau$ so that ${\bf A}_\tau(g_1,g_2) = \inp{g_1}{A_\tau g_2}_{L^2_\tau}$. Consequently, 
\begin{align}
    (A_{\tau}g)(t) 
    &= \frac{1}{\tau}\inp[\Big]{\int_{0}^{\tau}g(s)\mathcal{K}^sf\mathtt{d}s}{\mathcal{K}^tf}_{\mathcal{H}_f} = \frac{1}{\tau}\displaystyle\int_{0}^\tau g(s)\rho(t-s)\mathtt{d}s,
\end{align}
which is the same as the definition of $A_\tau$ in \cite{Vautard1989SingularSA}. Assumption \ref{assumption: K continuous} implies that $\rho\in L^2_{\text{loc}}(\mathbb{R})$. This implies that $A_{\tau}$ is a Hilbert-Schmidt operator on $L^2_\tau$.
We shall use the following variational description of the eigenvalues.
\begin{proposition}[The min-max principle]\label{proposition: minmax}
Let $\mathcal{H}$ be a Hilbert space and $A$ a Hermitian operator on $\mathcal{H}$. Let $\lambda_1\geq \lambda_2\geq \cdots$ be the eigenvalues of $A$ in descending order. Then
\begin{align}
    \lambda_i = \max_{\substack{\mathcal{M}\subset\mathcal{H}\\\dim\mathcal{M} = i}}\min_{v\in\mathcal{M}} \frac{\inp{v}{Av}}{\|v\|^2}
\end{align}
\end{proposition}

Our main result states that,
\begin{theorem}[Main result]\label{theorem: main result}
Under assumptions \ref{assumption: stationary} and \ref{assumption: K continuous}, we have, for all $i\in\mathbb{N}$
\begin{align}
    \lim_{\tau\to\infty}\lambda_{\tau,i} = |a_i|^2,
\end{align}
where $\lambda_{\tau,i}$ stands for the eigenvalues of $A_{\tau}$.
\end{theorem}

The following proposition \cite{Zhen2021} demonstrates the correspondence between the eigenfrequencies of the continuous-time time-shift operator and the discrete-time time-shift operator. Please refer to \cite{Zhen2021} for the notations in the proposition.
\begin{proposition}
 Let $\{f(X_t):t\geq 0\}$ be a continuous time process for which $\rho_s$ exists for all $s\geq 0$. Let $\Delta t>0$ be a time step. Assume that 
\begin{align}
    &\lim_{T\to\infty}\frac{1}{T}\int_{0}^Tf(X_t)\bar{f}(X_{t+k\Delta t})dt\nonumber \\
    =&\lim_{T\to\infty}\frac{\Delta t}{T}\sum_{\mathbb{N}\ni n=0}^{T/\Delta t}f(X_{n\Delta t})\bar{f}(X_{(n+k)\Delta t}),\label{eq: inner product cont=disc}
\end{align}
for all $k\in\mathbb{N}$. Then $\mathcal{H}_f \hookrightarrow \mathcal{H}^{\text{cont}}_f$.
Let $q$ be an eigenfrequency of the discrete-time operator $\mathcal{K}^{\Delta t}$, i.e. there exists $h\in \mathcal{H}_f \hookrightarrow \mathcal{H}^{\text{cont}}_f$ so that $\mathcal{K}^{\Delta t}h = e^{iq}h$. Then there exists an integer $k$, and $h_k\in\mathcal{H}^{\text{cont}}_f$, so that
\begin{align}
    \mathcal{K}^sh_k = e^{i\frac{q+2k\pi}{\Delta t}s}h_k
\end{align}
for all $s\geq 0$.
\end{proposition}

\begin{remark}
It is worth to point out that the one-parameter semigroup of isometries $\{\mathcal{K}^s: s\geq 0\}$ is equivalent to the classical Koopman one-parameter semigroup $\{\tilde{\mathcal{K}}^s: s\geq 0\}$ which acts on $L^2(X,d\nu)$ almost surely (with respect to the initial state of the time series),  if the dynamical system is ergodic and has finite invariant measure $\nu$ on the phase space $X$. Because if $f\in L^2(X,\nu)$, then $f\tilde{\mathcal{K}}^s\bar{f}\in L^1(X,d\nu)$ and Birkhoff ergodic theorem states that $\rho(s) = \nu(f\tilde{\mathcal{K}}^s\bar{f})$ for almost every initial state $x_0\in X$. In other words, $\inp{f}{\mathcal{K}^sf}_{\mathcal{H}_f} = \inp{f}{\tilde{\mathcal{K}}^sf}_{L^2(X,d\nu)}$. Note that $f$ is interpreted as a given time series on the left of the equality and interpreted as a function on the right of the equality. This shows that under the assumption that the dynamical system is ergodic and (finite) measure-preserving, there is a natural isometric bijection from $\mathcal{H}_f$ to $L^2(X, d\nu)$. 
\end{remark}
For mathematical interests, we present the main result in an abstract mathematical form.
\begin{theorem}[Main result in mathematical form]
Let $\mathcal{H}$ be a Hilbert space and $\{\mathcal{K}^s: s\geq 0\}$ a strongly continuous one-parameter semigroup of isometries acting on $\mathcal{H}$. For any $f\in\mathcal{H}$, let $f = \displaystyle\sum_{i}a_iv_i + f^{\perp}$, where $v_i$'s are the common eigenvectors of $\mathcal{K}^s$ for all $s\geq 0$, and $f^{\perp}$ is the component of $f$ that is orthogonal to the eigenspace of $\mathcal{K}^s$ for all $s\geq 0.$ Assume that $|a_1|\geq |a_2|\geq \cdots \geq 0$. For any $\tau > 0$, let $A_{f,\tau}$ be the Hermitian operator on $L^2([0,\tau])$, such that for any $g\in L^2([0,\tau])$ and any $t\in[0,\tau]$,
\begin{align}
    (A_{f,\tau}g)(t) =\frac{1}{\tau} \int_{0}^\tau g(s)\inp{\mathcal{K}^sf}{\mathcal{K}^tf}_{\mathcal{H}}\mathtt{d}s.
\end{align}
Then $A_{f,\tau}$ is a Hilbert-Schmidt operator and hence has purely punctual spectrum. Let $\lambda_{f,\tau,i}$ be the $i-$th largest eigenvalue of $A_{f,\tau}$. Then we have
\begin{align}
    \lim_{\tau\to\infty}\lambda_{f,\tau,i} = |a_i|^2.
\end{align}
\end{theorem}

\section{Proof of the main result}
For any fixed small $\epsilon\geq 0$, we choose $k$, so that $\displaystyle\sum_{i=k+1}^{\infty}|a_i|^2 \leq \epsilon$.
We have the orthogonal decomposition 
\begin{align}
    f =& f_d + f_{NU} + f_{c} = \sum_{i=1}^ka_iv_i + \sum_{i=k+1}^{\infty}a_iv_i + f_{d,k} + f_{NU} + f_c \nonumber \\
    =& f_{d,k} + f_{d,\epsilon} + f_{NU} + f_{c},
\end{align}
where $f_{d,k}\in\mathcal{H}_{f,d,k}$ which is the subspace of $\mathcal{H}_{f,d}$ spanned by $\{v_1,...,v_k\}$,
and $f_{d,\epsilon}\in \mathcal{H}_{f,d,\epsilon}$ the subspace spanned by the rest of the eigenvectors, $f_{NU}\in \mathcal{H}_{f,NU}$, and $f_{c}\in\mathcal{H}_{f,c}$. Note that $\mathcal{H}_{f,d,k}$, $\mathcal{H}_{f,d, \epsilon}$, $\mathcal{H}_{f,NU}$, and $\mathcal{H}_{f,c}$ are pairwise orthogonal and invariant subspaces of $\mathcal{H}_f$. Hence following Eq.\eqref{eq: <gAg> = <int,int>}, for any $g_1,g_2\in L^2_{\tau}$, 
\begin{align}
    &\inp{g_1}{A_{\tau}g_2}_{L^2_\tau} =  \frac{1}{\tau}\inp[\Big]{\int_{0}^{\tau}g_1(s)\mathcal{K}^sf\mathtt{d}s}{\int_{0}^{\tau}g_2(t)\mathcal{K}^tf\mathtt{d}t}_{\mathcal{H}_f}\nonumber \\
    =& \frac{1}{\tau} \inp[\Big]{\int_{0}^\tau g_1(s)\mathcal{K}^s(f_{d,k} + f_{d,\epsilon} + f_c + f_{NU})\mathtt{d}s}{\int_{0}^\tau g_2(t)\mathcal{K}^t(f_{d,k} + f_{d,\epsilon} + f_c + f_{NU})\mathtt{d}t}\nonumber \\
    =& \frac{1}{\tau}\inp[\Big]{\int_{0}^{\tau}g_1(s)\mathcal{K}^sf_{d,k}\mathtt{d}s}{\int_{0}^{\tau}g_2(t)\mathcal{K}^tf_{d,k}\mathtt{d}t}_{\mathcal{H}_f}\nonumber \\
    &+\frac{1}{\tau}\inp[\Big]{\int_{0}^{\tau}g_1(s)\mathcal{K}^sf_{d,\epsilon}\mathtt{d}s}{\int_{0}^{\tau}g_2(t)\mathcal{K}^tf_{d,\epsilon}\mathtt{d}t}_{\mathcal{H}_f}\nonumber \\
    &+\frac{1}{\tau} \inp[\Big]{\int_{0}^{\tau}g_1(s)\mathcal{K}^sf_c\mathtt{d}s}{\int_{0}^{\tau}g_2(t)\mathcal{K}^tf_c\mathtt{d}t}_{\mathcal{H}_f}\nonumber \\
    &+ \frac{1}{\tau}\inp[\Big]{\int_{0}^{\tau}g_1(s)\mathcal{K}^sf_{NU}\mathtt{d}s}{\int_{0}^{\tau}g_2(t)\mathcal{K}^tf_{NU}\mathtt{d}t}_{\mathcal{H}_f}\nonumber \\
    =&\inp{g_1}{A_{\tau,d,k}g_2}_{L^2_\tau} +\inp{g_1}{A_{\tau,d,\epsilon}g_2}_{L^2_\tau} + \inp{g_1}{A_{\tau,c}g_2}_{L^2_\tau} + \inp{g_1}{A_{\tau,NU}g_2}_{L^2_{\tau}},
\end{align}
in which the definition of $A_{\tau,d,k}$, $A_{\tau,d,\epsilon}$, $A_{\tau,c}$ and $A_{\tau,NU}$ are obvious. It is not hard to show that $A_{\tau,d,k}$, $A_{\tau,d,\epsilon}$, $A_{\tau,c}$ and $A_{\tau,NU}$ all admit eigendecomposition since they are all Hilbert-Schmidt Hermitian operators. Note that the  cross product terms all as  $\mathcal{H}_{f,d,k}$, $\mathcal{H}_{f,d,\epsilon}$, $\mathcal{H}_{f,c}$ and $\mathcal{H}_{f,NU}$ are pairwise orthogonal and invariant under $\mathcal{K}^s$ for all $s\geq 0$.

Let $\lambda_{\tau_d,k,i}$, $\lambda_{\tau,d,\epsilon,i}$, $\lambda_{\tau,c,i}$, and $\lambda_{\tau,NU,i}$ be the $i-$th largest eigenvalue of $A_{\tau,d,k}$, $A_{\tau,d,\epsilon}$, $A_{\tau,c}$, $A_{\tau,NU}$ respectively. We will prove the following identities:
\begin{proposition}\label{proposition: lim lambda for dk,de,c,NU}
\begin{align}
    &\lim_{\tau\to\infty}\lambda_{\tau,d,k,i} = |a_i|^2 \text{ for $i=1,...,k$,}\label{eq: lim lambda dk}\\
    &\lambda_{\tau,d,\epsilon,1} \leq \epsilon \text{ for any $\tau>0$,} \label{eq: lim lambda de}\\
    &\lim_{\tau\to\infty}\lambda_{\tau,c,1} = 0 \label{eq: lim lambda c},\\
    &\lim_{\tau\to\infty}\lambda_{\tau,NU,1} = 0 \label{eq: lim lambda NU}.
\end{align}
\end{proposition}
Before we start to prove Eqs.\eqref{eq: lim lambda dk}-\eqref{eq: lim lambda NU}, it is not hard to see that proposition \ref{proposition: minmax} and proposition \ref{proposition: lim lambda for dk,de,c,NU} directly implies the main result. Indeed, for any fixed $n$ and any $\epsilon>0$, we can find $k$ so that $n\leq k$ and $\displaystyle\sum_{i=k+1}^{\infty}|a_i|^2\leq \epsilon$. Then we find $\tau$ large enough so that $\lambda_{\tau,c,1}\leq \epsilon$ and $\lambda_{\tau,NU,1}\leq \epsilon$. Note that $A_{\tau,d,k}$, $A_{\tau,d,\epsilon}$, $A_{\tau,c}$, and $A_{\tau,NU}$ are all positive semi-definite. Applying the min-max principle we have
\begin{align}
    \lambda_{\tau,n} &= \max_{\substack{\mathcal{M}\subset L^2_{\tau} \\\dim\mathcal{M} = n}}\min_{v\in\mathcal{M}} \frac{\inp{v}{A_{\tau}\;v}}{\|v\|^2}\\
    &= \max_{\substack{\mathcal{M}\subset L^2_\tau\\\dim\mathcal{M} = n}}\min_{v\in\mathcal{M}}\frac{\inp{v}{A_{\tau,d,k}\;v} + \inp{v}{A_{\tau,d,\epsilon}\;v} + \inp{v}{A_{\tau,c}\;v} + \inp{v}{A_{\tau,NU}\;v}}{\|v\|^2}\\
    &\geq \max_{\substack{\mathcal{M}\subset L^2_\tau\\\dim\mathcal{M} = n}}\min_{v\in\mathcal{M}}\frac{\inp{v}{A_{\tau,d,k}\;v}}{\|v\|^2} = \lambda_{\tau,d,k,n},
\end{align}
and that 
\begin{align}
    \lambda_{\tau,n} &= \max_{\substack{\mathcal{M}\subset L^2_{\tau} \\\dim\mathcal{M} = n}}\min_{v\in\mathcal{M}} \frac{\inp{v}{A_{\tau}\;v}}{\|v\|^2}\\
    &\leq \max_{\substack{\mathcal{M}\subset L^2_{\tau} \\\dim\mathcal{M} = n}}\min_{v\in\mathcal{M}} \frac{\inp{v}{A_{\tau,d,k}\;v}}{\|v\|^2} + 2\epsilon = \lambda_{\tau,d,k,n}+2\epsilon.
\end{align}
Combined with Eq.\eqref{eq: lim lambda dk}, this implies Theorem \ref{theorem: main result}.
\begin{proof}[Eq.\eqref{eq: lim lambda dk}]
 Recall from Eq.\eqref{eq: phi(vi) = dirac} that each eigenvector $v_i$ corresponds to a point $x_i$ in the support of $\mu_{d}$. For any $g\in L^2_\tau$,
 Theorem \ref{theorem: spectral representation continuous} states that $\displaystyle\int_{0}^{\tau}g(s)\mathcal{K}^sf_{d,k}\mathtt{d}s$ has the following representation in $L^2(\mathbb{R},d\mu)$, for any $x\in\mathbb{R}$,
 \begin{align}
     \Big{(}\phi{\big (}\int_{0}^{\tau}g(s)\mathcal{K}^sf_{d,k}\mathtt{d}s\big{)}\Big{)}(x) = \begin{cases} \int_{0}^{\tau}g(s)e^{isx_j}\mathtt{d}s &\text{, if $x=x_j$ for some $j$.}\\
     0 &\text{, otherwise.}
     \end{cases}
 \end{align}
And
\begin{align}
    \inp{g}{A_{\tau,d,k}g}_{L^2_{\tau}} &= \frac{1}{\tau}\inp[\Big]{\int_{0}^{\tau}g(s)\mathcal{K}^sf_{d,k}\mathtt{d}s}{\int_{0}^{\tau}g(t)\mathcal{K}^tf_{d,k}\mathtt{d}t}_{\mathcal{H}_f}\\
    &= \frac{1}{\tau}\sum_{j=1}^k{\Big \|}\int_{0}^{\tau}g(s)e^{isx_j}\mathtt{d}s{\Big \|}^2_{L^2(\mathbb{R},d\mu)}\\
    &= \frac{1}{\tau}\sum_{j=1}^k|a_j|^2\Big{|}\int_{0}^{\tau}g(s)e^{isx_j}\mathtt{d}s\Big{|}^2.
\end{align}
Let $\xi_j\in L^2_{\tau}$ so that $\xi_j(s) = e^{isx_j}$ for any $s\in [0,\tau]$. Then $\|\xi_j\|^2_{L^2_\tau} = \tau$ and 
\begin{align}
    \inp{g}{A_{\tau,d,k}\;g}_{L^2_{\tau}} &= \frac{1}{\tau}\sum_{j=1}^k|a_j|^2|\inp{\xi_j}{g}_{L^2_\tau}|^2 = \sum_{j=1}^k|\inp{\frac{a_j\xi_j}{\sqrt{\tau}}}{g}_{L^2_\tau}|^2\label{eq: gAg = sum <axi g>**2}
\end{align}
Let $V_{\tau,k} = \text{Span}_{\mathbb{C}}\{\frac{a_1\xi_1}{\sqrt{\tau}}, \frac{a_2\xi_2}{\sqrt{\tau}},\cdots,\frac{a_k\xi_k}{\sqrt{\tau}}\}$. We write $g = g_{\tau,k} + g^{\perp}$, where $g_{\tau,k}\in V_{\tau,k}$, and $g^{\perp}\in V_{\tau,k}^{\perp}$. 
Then
\begin{align}
    \inp{g}{A_{\tau,k,d}\;g}_{L^2_\tau} = \sum_{j=1}^k|\inp{\frac{a_j\xi_j}{\sqrt{\tau}}}{g_{\tau,k}}|_{L^2_\tau}^2.
\end{align}
Note that $\dim V_{\tau,k} = k$ for all $\tau > 0$. Direct calculation yields that, for $j\neq \ell$, $\inp{\frac{a_j\xi_j}{\sqrt{\tau}}}{\frac{a_\ell\xi_\ell}{\sqrt{\tau}}}_{L^2_\tau} = a_j\bar{a}_l \frac{e^{i(x_j-x_\ell)\tau}-1}{i\tau(x_j-x_\ell)}\to 0$ as $\tau\to\infty$. Therefore the distribution of the eigenvalues of $A_{\tau,k,d}$ shall approach to the distribution of the eigenvalues of
\begin{align}
    \begin{pmatrix}
    |a_1|^2 & 0 & \cdots & 0 \\
    0 & |a_2|^2 & \cdots & 0\\
    \vdots\\
    0 & 0 & \cdots & |a_k|^2
    \end{pmatrix}
\end{align}
as $\tau\to\infty$. This completes the proof of Eq.\eqref{eq: lim lambda dk}.
\end{proof}

\begin{proof}[Eq.\eqref{eq: lim lambda de}]
Similar to Eq.\eqref{eq: gAg = sum <axi g>**2}, for any $g\in L^2_{\tau}$, $\|g\|_{L^2_\tau} = 1$, we have 
\begin{align}
    \inp{g}{A_{\tau,d,\epsilon}\;g}_{L^2_{\tau}} &= \frac{1}{\tau}\sum_{j=k+1}^\infty|a_j|^2|\inp{\xi_j}{g}_{L^2_\tau}|^2 = \sum_{j=k+1}^\infty|\inp{\frac{a_j\xi_j}{\sqrt{\tau}}}{g}_{L^2_\tau}|^2\\
    &\leq \sum_{j=k+1}^\infty |a_j|^2 \leq \epsilon.
\end{align}
Then the min-max principle implies that $\lambda_{\tau,k,\epsilon,1} \leq \epsilon$.
\end{proof}

\begin{proof}[Eq.\eqref{eq: lim lambda c}]
Following \cite{Halmos1956} (page 39-41), we first show that 
\begin{align}
    \lim_{\tau\to\infty}\frac{1}{\tau}\int_{0}^{\tau}\big{|}\mu_{f,c}(e^{isx})\big{|}\mathtt{d}s = 0, \label{eq: mean moments tends to 0 for continuous spectrum}
\end{align}
or equivalently
\begin{align}
    \lim_{\tau\to\infty}\frac{1}{\tau}\int_{0}^{\tau}\big{|}\mu_{f,c}(e^{isx})\big{|}^2\mathtt{d}s = 0.\label{eq: mean moments^2 tends to 0 for continuous spectrum}
\end{align}
Eq.\eqref{eq: mean moments tends to 0 for continuous spectrum} means that the large moments associated to the continuous spectral measure has density zero. For any $\epsilon>0$, we write 
$\mu_{f,c} = \mu_{f,c,1} + \mu_{f,c,\epsilon}$, in which $\mu_{f,c,1}$ has compact support, $\mu_{f,c,\epsilon}(\mathbb{R})<\epsilon$ and $\mu_{f,c,1}\perp\mu_{f,c,\epsilon}$. Denote the support of $\mu_{f,c,1}$ by $B_1$. Then we have
\begin{align}
    \frac{1}{\tau}\int_{0}^\tau \big{|}\mu_{f,c}(e^{isx})\big{|}^2\mathtt{ds} &= \frac{1}{\tau}\int_{0}^\tau \big{|}\mu_{f,c,1}(e^{isx})\big{|}^2\mathtt{ds} + \frac{1}{\tau}\int_{0}^\tau \big{|}\mu_{f,c,\epsilon}(e^{isx})\big{|}^2\mathtt{ds}\\
    &< \frac{1}{\tau} \int_{0}^\tau\big{|}\int_{\mathbb{R}}e^{isx}\mathtt{d}\mu_{f,c,1}(x)\big{|}^2\mathtt{d}s + \epsilon
\end{align}
and that 
\begin{align}
    &\frac{1}{\tau}\int_{0}^\tau \big{|}\mu_{f,c,1}(e^{isx})\big{|}^2\mathtt{ds} = \frac{1}{\tau} \int_{0}^\tau\big{|}\int_{\mathbb{R}}e^{isx}\mathtt{d}\mu_{f,c,1}(x)\big{|}^2\mathtt{d}s \nonumber \\
    =&\frac{1}{\tau}\int_{0}^\tau\mathtt{d}s\int_{\mathbb{R}}\int_{\mathbb{R}}e^{is(x-y)}\mathtt{d}\mu_{f,c,1}(x)\mathtt{d}\mu_{f,c,1}(y)\\
    =&\frac{1}{\tau}\int_{\mathbb{R}}\int_{\mathbb{R}}\mathtt{d}\mu_{f,c,1}(x)\mathtt{d}\mu_{f,c,1}(y)\int_{0}^\tau e^{is(x-y)}\mathtt{d}s\\
    =&\frac{1}{\tau}\int_{B_1}\int_{B_1}\mathtt{d}\mu_{f,c,1}(x)\mathtt{d}\mu_{f,c,1}(y)\int_{0}^\tau e^{is(x-y)}\mathtt{d}s\label{int-prod}
\end{align}
Note that $\Big{|}\frac{1}{\tau}\int_{0}^{\tau}e^{is(x-y)}\mathtt{d}s\Big{|} \leq 1$
for any $\tau>0$ and any $x,y\in\mathbb{R}$. And when $x\neq y$ 
\begin{align}
\label{limit}
 1\geq \Big{|}\frac{1}{\tau}\int_{0}^{\tau}e^{is(x-y)}\mathtt{d}s\Big{|} =\Big{|}\frac{e^{i\tau (x-y)}-1}{\tau i(x-y)}\Big{|}\underset{\tau\to\infty}{\longrightarrow} 0. 
\end{align}
Since $\mu_{f,c,1}$ is continuous, we have  that $(\mu_{f,c,1}\times\mu_{f,c,1})(\{(x,y)\in\mathbb{R}^2: x=y\}) = 0$. Hence, the integral in Eq.\eqref{int-prod} boils down to an integral on $\mathbb{R}^2\setminus\{x=y\}$. Lebesgue's dominated convergence theorem implies that the integral in Eq.\eqref{int-prod} converges to 0 as $\tau\to\infty$. Hence $\displaystyle\limsup_{\tau\to\infty}\frac{1}{\tau}\int_{0}^\tau\big{|}\mu_{f,c}(e^{isx})\big{|}^2\mathtt{d}s <\epsilon$ for any $\epsilon>0$. This implies Eq.\eqref{eq: mean moments^2 tends to 0 for continuous spectrum}.

For any $g\in L^2(\mathbb{R})$, Theorem \ref{theorem: spectral representation continuous} implies that 
\begin{align}
    \Big{(}\phi\big{(}\int_{0}^\tau g(s)\mathcal{K}^sf_{d,c}\mathtt{d}s\big{)}\Big{)}(x) = \int_{0}^{\tau}g(s)e^{isx}\mathtt{d}s.
\end{align}
Therefore 
\begin{align}
    \inp{g}{A_{\tau,c}\;g}_{L^2_\tau} &= \frac{1}{\tau}\inp[\Big]{\int_{0}^{\tau}g(s)\mathcal{K}^sf_{d,c}\mathtt{d}s}{\int_{0}^{\tau}g(t)\mathcal{K}^tf_{d,c}\mathtt{d}t}_{\mathcal{H}_f}\\
    &= \frac{1}{\tau}\inp[\Big]{\phi\big{(}\int_{0}^{\tau}g(s)\mathcal{K}^sf_{d,c}\mathtt{d}s\big{)}}{\phi\big{(}\int_{0}^{\tau}g(t)\mathcal{K}^tf_{d,c}\mathtt{d}t\big{)}}_{L^2(\mathbb{R},d\mu_c)}\\
    &= \frac{1}{\tau} \int_{-\infty}^{\infty}\mathtt{d}\mu_{f,c}(x)\int_{0}^\tau\int_{0}^\tau g(s)\bar{g}(t)e^{i(s-t)x}\mathtt{d}s\mathtt{d}t\\
    &=\frac{1}{\tau}\int_{0}^\tau\int_{0}^\tau g(s)\bar{g}(t)\mu_{f,c}\;(e^{i(s-t)x})\mathtt{d}s\mathtt{d}t
\end{align}
Hence
\begin{align}
    |\inp{g}{A_{\tau,c}\;g}|\leq& \frac{1}{\tau}\int_{0}^\tau\int_0^\tau |g(t)|\cdot|g(s)|\cdot |\mu_{f,c}\;(e^{i(s-t)x})|\mathtt{d}t\mathtt{d}s\\
    =&\frac{1}{\tau}\iint_{0\leq s\leq t\leq \tau} |g(t)|\cdot|g(s)|\cdot |\mu_{f,c}\;(e^{i(s-t)x})|\mathtt{d}t\mathtt{d}s\\
    &+\frac{1}{\tau}\iint_{0\leq t\leq s\leq \tau} |g(t)|\cdot|g(s)|\cdot |\mu_{f,c}\;(e^{i(s-t)x})|\mathtt{d}t\mathtt{d}s\\
    =&\frac{2}{\tau}\int_{0}^\tau |g(t)|\int_{t}^\tau |g(s)|\cdot|\mu_{f,c}\;(e^{i(s-t)x})|\mathtt{d}t\mathtt{d}s\\
    =& \frac{2}{\tau}\int_{0}^{\tau}|g(t)|\int_{0}^{\tau-t}|g(t+s)|\cdot|\mu_{f,c}(e^{isx})|\mathtt{d}s\mathtt{dt}\\
    \leq&\frac{2}{\tau} \int_{0}^{\tau}\int_{0}^{\tau-t}\frac{1}{2}(|g(t)|^2+|g(t+s)|^2)|\mu_{f,c}(e^{isx})|\mathtt{d}s\mathtt{d}t\\
    =&\frac{1}{\tau}\int_{0}^{\tau}|\mu_{f,c}(e^{isx})|\int_{0}^{\tau-s}(|g(t)|^2+|g(t+s)|^2)\mathtt{d}s\mathtt{d}t\\
    \leq& \frac{1}{\tau}\int_0^\tau 2|\mu_{f,c}(e^{isx})|\cdot\|g\|^2_{L^2_{\tau}}\mathtt{d}s
\end{align}
Therefore 
\begin{align}
    \lambda_{\tau,c,1} = \max_{g\in L^2_\tau}\frac{\inp{g}{A_{\tau,c}g}}{\|g\|_{L^2_\tau}}\to 0,
\end{align}
as $\tau\to\infty$. This completes the proof of Eq.\eqref{eq: lim lambda c}.
\end{proof}

\begin{proof}[Eq.\eqref{eq: lim lambda NU}]
Recall that $\mathcal{H}_{f,NU} \cong L^2([0,+\infty], \mathcal{H}_0)$. Hence $f_{NU}$ can be represented as a curve  from $[0,\infty]$ to $\mathcal{H}_0$. We denote this curve by $\gamma$. Without ambiguity, we do not distinguish between $\gamma$ and $f_{NU}$. Hence for each $t\geq 0$, $\gamma(t)\in\mathcal{H}_0$. And $\|\gamma\|^2_{\mathcal{H}_{f,NU}} = \displaystyle\int_{0}^\infty\|\gamma(t)\|_{\mathcal{H}_0}^2\mathtt{d}t$. Recall that $(\mathcal{K}^s\gamma)(t) = \gamma(t+s)$. We set $\gamma(t)=0$ for all $t<0$. Hence for any $g\in L^2_\tau$,
\begin{align}
    &\inp{g}{A_{\tau,NU}\;g}_{L^2_\tau} = \frac{1}{\tau}\inp[\Big]{\int_{0}^{\tau}g(s_1)\mathcal{K}^{s_1}\gamma\mathtt{d}s_1}{\int_{0}^{\tau}g(s_2)\mathcal{K}^{s_2}\gamma\mathtt{d}s_2}_{\mathcal{H}_{f,NU}}\\
    =&\frac{1}{\tau}\int_{0}^{\infty}\int_{0}^\tau\int_{0}^\tau\bar{g}(s_2)g(s_1)\inp{\gamma(t+s_1)}{\gamma(t+s_2)}_{\mathcal{H}_0}\mathtt{d}s_1\mathtt{d}s_2\mathtt{d}t\\
    =& \frac{1}{\tau}\int_{0}^\tau\int_{0}^\tau\bar{g}(s_2)g(s_1)\int_{0}^{\infty}\inp{\gamma(t+s_1)}{\gamma(t+s_2)}_{\mathcal{H}_0}\mathtt{d}t\mathtt{d}s_1\mathtt{d}s_2\label{eq: tmp 1}
\end{align}
We first show the following identity:
\begin{align}
    \lim_{s\to\infty}\inp{\gamma}{\mathcal{K}^s\gamma}_{\mathcal{H}_{f,NU}} = \lim_{s\to\infty}\int_{0}^\infty\inp{\gamma(t)}{\gamma(t+s)}_{\mathcal{H}_0}\mathtt{d}t = 0.\label{eq: NU lim int <gamma Ksgamma> = 0}
\end{align}
To prove Eq.\eqref{eq: NU lim int <gamma Ksgamma> = 0}, without loss of generality we assume that $\|\gamma\|_{\mathcal{H}_{f,NU}} = 1$. For any $\epsilon > 0$, there exists $N_{\epsilon}$, so that $\int_{0}^{N_{\epsilon}}\|\gamma(t)\|_{\mathcal{H}_0}^2\mathtt{d}t > 1-\epsilon$. This means that $\int_{N_{\epsilon}}^{\infty}\|\gamma(t)\|^2\mathtt{d}t < \epsilon$. Therefore for any $s\geq N_{\epsilon}$,
\begin{align}
    \Big{|}\int_{0}^{\infty}\inp{\gamma(t)}{\gamma(t+s)}_{\mathcal{H}_0}\mathtt{d}t\Big{|}^2\leq \Big{|}\int_{0}^{\infty}\|\gamma(t)\|^2_{\mathcal{H}_0}\mathtt{d}t\Big{|}^2\cdot\Big{|}\int_{N_{\epsilon}}^{\infty}\|\gamma(t)\|^2_{\mathcal{H}_0}\mathtt{d}t\Big{|}^2<\epsilon^2.
\end{align}
This proves Eq.\eqref{eq: NU lim int <gamma Ksgamma> = 0}. Now we continue with Eq.\eqref{eq: tmp 1}:
\begin{align}
    &\inp{g}{A_{\tau,NU}g}_{L^2_\tau} \leq \Big{|}\frac{2}{\tau}\int_{0}^{\tau}\int_{s_1}^{\tau}\bar{g}(s_2)g(s_1)\inp{\mathcal{K}^{s_1}\gamma}{\mathcal{K}^{s_2}\gamma}_{\mathcal{H}_{f,NU}}\mathtt{d}s_1\mathtt{d}s_2\Big{|}\\
    \leq&\Big{|}\frac{2}{\tau}\int_{0}^{\tau}\int_{0}^{\tau-s_1}g(s_1)\bar{g}(s_1+s)\inp{\gamma}{\mathcal{K}^{s}\gamma}_{\mathcal{H}_{f,NU}}\mathtt{d}s_1\mathtt{d}s\Big{|}
\end{align}
For any $\epsilon>0$, find $M_{\epsilon}$, so that for any $|\inp{\gamma}{\mathcal{K}^s\gamma}| < \epsilon$ for any $s> M_{\epsilon}$. Now for any $\tau>M_{\epsilon}/\epsilon$ and any $\|g\|_{L^2_\tau} = 1$, we have
\begin{align}
    &\inp{g}{A_{\tau,NU}g}_{L^2_\tau}\\ \leq&\frac{2}{\tau}\int_{0}^\tau\int_{0}^{M_{\epsilon}}|g(s_1)|\cdot |g(s_1+s)|\cdot|\inp{\gamma}{\mathcal{K}^{s}\gamma}_{\mathcal{H}_{f,NU}}|\mathtt{d}s_1\mathtt{d}s + \\
    &\frac{2}{\tau}\int_{0}^\tau\int_{M_{\epsilon}}^{\tau-s_1}|g(s_1)|\cdot |g(s_1+s)|\cdot|\inp{\gamma}{\mathcal{K}^{s}\gamma}_{\mathcal{H}_{f,NU}}|\mathtt{d}s_1\mathtt{d}s\\
    \leq&\frac{1}{\tau}\int_{0}^{\tau} \int_{0}^{M_{\epsilon}}(|g(s_1)|^2 + |g(s_1+s)|^2)|\inp{\gamma}{\mathcal{K}^{s}\gamma}_{\mathcal{H}_{f,NU}}|\mathtt{d}s_1\mathtt{d}s + \\
    &\frac{1}{\tau}\int_{0}^{\tau} \int_{M_{\epsilon}}^{\tau-s_1}(|g(s_1)|^2 + |g(s_1+s)|^2)|\inp{\gamma}{\mathcal{K}^{s}\gamma}_{\mathcal{H}_{f,NU}}|\mathtt{d}s_1\mathtt{d}s \\
    \leq& \frac{1}{\tau}\int_0^\tau\int_0^{M_{\epsilon}}|g(s_1)|^2\mathtt{d}s_1\mathtt{d}s + \frac{1}{\tau}\int_0^\tau\int_0^{M_{\epsilon}}|g(s_1+s)|^2\cdot|\inp{\gamma}{\mathcal{K}^s\gamma}_{\mathcal{H}_{f,NU}}|\\
    &\mathtt{d}s_1\mathtt{d}s+\frac{1}{\tau}\int_{0}^{\tau} \int_{M_{\epsilon}}^{\tau-s_1}\epsilon(|g(s_1)|^2 + |g(s_1+s)|^2)\mathtt{d}s_1\mathtt{d}s\\
    \leq& \frac{M_{\epsilon}}{\tau} + \frac{M_{\epsilon}}{\tau} + 2\frac{\epsilon}{\tau}(\tau-M_\epsilon) \leq 4\epsilon.
\end{align}
Therefore for $\tau>M_{\epsilon}/\epsilon$,
\begin{align}
    \lambda_{\tau,NU,1} = \max_{\substack{g\in L^2_\tau\\\|g\|=1}}\inp{g}{A_{\tau,NU}g}\leq 4\epsilon.
\end{align}
This completes the proof of Eq.\eqref{eq: lim lambda NU}.
\end{proof}

%
%
%
\providecommand{\noopsort}[1]{}\providecommand{\singleletter}[1]{#1}%

\end{document}